# Enhancing Flexibility in Smart Manufacturing: System Property Enabled Multiagent Approach for Mobile Robot Scheduling in Multiproduct Flexible Manufacturing Systems

Muhammad Waseem[1], Qing Chang[1*]

*1 Department of Mechanical and Aerospace Engineering, University of Virginia, Charlottesville VA 22903*

*Corresponding author, Email: qc9nq@virginia.edu

**Abstract**

The present volatile global market is marked by a rising need for extensively tailored products. Flourishing in such an environment necessitates fostering a more intimate connection between market requirements and the manufacturing system, with a primary emphasis on the customer. Multiproduct flexible manufacturing systems (FMS) appear to be one of the most suitable manufacturing systems for mass-customization, given their capacity to readily adapt to variations in product specifications and functionalities to meet evolving market demands. Mobile robots play a pivotal role in managing multiproduct FMS. However, scheduling these robots poses a challenging problem due to system internal uncertainties and fluctuations in external demand. This paper proposes a solution to the robots' scheduling problem in a multiproduct FMS with stochastic machine failures and uncertain market demand for different product types. The proposed approach abstracts system properties for FMS and leverages these properties to improve the multiagent deep deterministic policy gradient (MADDPG) algorithm. This modification involves integrating parallel layers into the critic network, focusing specifically on crucial information such as varying customer demands. To validate the effectiveness of the proposed method, a comprehensive comparative analysis is conducted, comparing it with the simple-MADDPG, DDPG, and DQN algorithms. The results of the analysis demonstrate an average improvement of 9% in training time and 19% in market demand satisfaction. Moreover, to ensure scalability, transfer learning is employed, showcasing the applicability of the proposed method in large-scale complex environments.

1. Introduction

The growing demand for personalized, intelligent, and sustainable products has driven the swift expansion of industrial Internet and cyber-physical technologies, leading to a transformative impact on manufacturing systems [1]. The industrial landscape has shifted from conventional manufacturing systems to smart manufacturing systems, utilizing self-contained, automated, and intelligent modules to harness flexibility. The growing emphasis on enhancing production planning and control is a response to evolving customer preferences for top-quality finished goods at competitive prices and shorter time-to-market. This shift motivates factories to excel and surpass competitors by adopting flexible manufacturing systems (FMS) [2]. The FMS has gained significant recognition in recent years as a pivotal concept characterized by its responsiveness and a heightened ability to adapt promptly and effectively [3]. The escalating market demand for customized products stands as a decisive catalyst driving shifts in manufacturing paradigms [4]. To meet rapidly evolving customization needs and uphold their competitive edge, manufacturers consistently enhance the flexibility of their production processes.

FMS, being a smart manufacturing system, requires robust control. Based on the literature, two major approaches are used in FMS control: centralized and decentralized control. Centralized control involves a single component overseeing all manufacturing entities, while decentralized control systems are tailored assemblies of autonomous components capable of intercommunication [2]. Decentralized control holds the potential to create a more efficient production environment, particularly when faced with deviations from standard operating conditions[5]. Consequently, there is a rising interest in flexible manufacturing systems and solutions that facilitate autonomous communication and operation of independent entities, such as multi-agent systems (MASs). [6]. However, the scheduling problem in FMS being of an NP-Hard nature is always a major challenge. Even though it has been addressed by several studies [7], its complexity increases when accompanied by other issues like stochasticity resulting from random machine disruptions, processing multiple product types and uncertain market demands. Besides, the complexity also depends on the type of FMS. While there have been several solutions proposed to address the common FMS like serial production lines, a mobile robot-assisted multiproduct FMS, where each product type follows a unique processing path, requires an intelligent control mechanism to deal with the system uncertainties.

Current approaches to FMS scheduling primarily fall into two categories [8]: exact methods and approximate methods. Exact methods, known for their robust optimization capabilities, are highly effective for tackling small-scale problems like mixed integer linear programming but come with a trade-off of high computational complexity. Approximate methods, including heuristic methods[9], meta-heuristic methods [10], and learning-based methods [11], offer satisfactory solutions within a shorter timeframe. Heuristic methods, exemplified by priority dispatching rules (PDRs), demonstrate rapid computation but may lack foresight due to limited consideration of problem properties. Meta-heuristic methods, drawing inspiration from natural phenomena, enhance solution quality through continuous iterative search from initial solutions. However, the complexity of designing meta-heuristic methods and the unpredictability of solution stability, influenced by random factors in the algorithm, pose challenges. Learning-based methods [11], categorized as supervised learning and reinforcement learning (RL), face complexities in obtaining data labels for the former, as they require optimal or approximate optimal solutions as labels. In contrast, the latter, especially deep reinforcement learning (DRL), emerges as a research focal point for FMS scheduling due to its rapid solving speed and robust generalization capabilities.

While the DRL algorithms can address complex problems and learn from extensive information, they may not be able to enjoy the advantage of the system property or expert knowledge about the system, resulting in prolonged training periods. Furthermore, these algorithms treat all information in the state space equally, making it difficult for the agent to capture crucial state information in the state space, explore solution space effectively in the definition of action space, and quickly identify key action in the design of reward function. It could be time-consuming to train an RL-based control policy considering all parameters and input data, assuming no knowledge about the system. Therefore, this paper studies a DRL multiagent deep deterministic policy gradient (MADDPG) algorithm and enhances the network structure to discern information importance levels within the state space, prioritizing critical details to expedite the training process and enhance decision-making efficiency.

The main contributions in this paper can be summarized as follows:

- We formulate mobile robot scheduling as a multiagent control problem to adapt dynamically to fluctuations in market demand, ensuring efficient production planning.

- We develop a novel MADDPG algorithm, incorporating the system property of opportunity window and market demand [12] through parallel layers in the critic network. It allows the model to capture complex relationships and dependencies between different products in a multiagent multiproduct FMS and improves the learning and decision-making capabilities.

- We employ transfer learning to enable the scalability of the proposed control strategies and demonstrate their applicability to large-scale environments.

The rest of the paper is organized as follows: Section 2 provides a literature review, and Section 3 delves into the system description. Formulation of the control problem is presented in Section 4. Section 5 covers MADDPG and its modification, while Section 6 introduces a case study to validate the proposed method. Results are discussed in Section 7. Section 8 explores transfer learning to assess the applicability of the proposed method in complex environments. Finally, Section 9 concludes the paper and outlines future work.

2. Literature review

There is extensive literature available on FMS, covering its various aspects such as modelling, analysis, and control [13]. The degree of flexibility an FMS can offer is inherently tied to its design and characteristics, with scheduling playing a pivotal role in control. Recent research has explored diverse methodologies like Petri-net [14], reinforcement learning (RL) [15], and deep reinforcement learning (DRL) [16] to address FMS management. For instance, Hu et al. [14] combined petri-net, DRL, and graph convolution networks (GCN) in a study focused on a dynamic FMS with shared resources. However, this approach may have limitations when applied to different environments. In another study, Windmann et al. [17] utilized a dynamic programming-based model to address routing issues in FMS, attempting to generalize the model for variable parametric settings. Despite these efforts, challenges persist, particularly in adapting these models to diverse environments.

Traditionally, scheduling literature has predominantly focused on jobs or machine scheduling, such as parallel machines [18], single machines [19], and job shops [20]. Conventional machine scheduling models assume the constant availability of material handling systems, neglecting material handling or transportation times during scheduling [21]. This assumption, however, does not hold in multiproduct FMS, where production output is intricately linked with robot scheduling and system modelling [22]. The dynamic nature of the system poses challenges for closed-form representations, particularly in the context of efficient scheduling in multiproduct FMS. The complexity arises from diverse product types, each with unique processing flows, random disruptions, uncertain demand fluctuations, and a multitude of machines, buffers, and robots.

In general, manufacturing systems have received significant attention in the literature [23], particularly for serial production lines [24]. However, it's essential to note that a multiproduct FMS differs significantly because each product type follows a distinct processing path. As customization trends surge, enterprises must promptly adapt production plans to meet evolving external demand, enhancing overall customer satisfaction [25, 26]. Our previous research [26] focused on modelling and controlling a multiproduct FMS, considering random machine disruptions and assuming fixed market demand for each product type. Yet, since demand is not always deterministic and can vary, accurate estimation becomes crucial as it significantly impacts customer satisfaction. While mathematical optimization approaches like metaheuristics and evolving models are traditionally employed in such scenarios [27], their efficacy diminishes due to the stochastic and NP-hard nature of the problem. The intricate relationship between changing production plans, different product paths, and fluctuating demand compounds the complexity of scheduling in multiproduct FMS. Relying solely on these approaches might not adequately address these multifaceted challenges.

Alternatively, machine learning-based approaches like Reinforcement Learning (RL) and Deep Reinforcement Learning (DRL) offer promising solutions for handling scheduling complexities in multiproduct FMS [12]. They've demonstrated success in learning optimal control policies in dynamic environments by iteratively learning from interactions. Their adaptability aligns with the needs of multiproduct FMS scheduling, yet their applicability is limited in contexts where multiple robots require coordination and interaction due to their focus on single-agent systems [28]. To address this, utilizing multiagent systems becomes crucial for accurate representation of real-world scenarios, enabling effective collaboration and information exchange among robots to enhance scheduling and decision-making [29]. Similarly, traditional policy gradient (PG) algorithms [30] are unsuitable for multi-agent environments due to distinct strategies among agents causing constant fluctuations, rendering convergence challenging when trained solely from a single-agent perspective. Thus, multiagent algorithms like multiagent actor-critic (MAAC) and MADDPG are used to overcome these limitations. MADDPG, with its specialized ability to handle interactions effectively, stands out by enabling cooperative behavior among agents in complex environments. [31].

MADDPG has been extensively applied across various domains, including autonomous vehicle control [32], traffic control [33], multi-agent games [34], and managing autonomous mobile robots (AMRs) [35]. However, depending on the problem's complexity, it may require longer training periods. In response, several derivative works have emerged, aiming to modify the MADDPG algorithm to suit specific problems. For instance, [36] introduced an actor-hierarchical attention critic to improve agent cooperation and facilitate efficient learning in mixed tasks. Additionally, [37] proposed a multiagent hierarchical DDPG to merge the advantages of multi-robot and hierarchical systems in DRL. Their primary objective involved decomposing complex tasks into simpler ones for corresponding solutions. In a cooperative and complex environment, [38] proposed the parameter-sharing DPG method, yet this method leads to increased computational complexity. Meanwhile, [39] introduced a decentralized generative policy network designed for effective agent cooperation in partially observable environments. Their approach involved incorporating an additional actor network within each agent's structure to replicate other agents' actions, eliminating the need to merge agents' target policies during training. Additionally, they introduced global joint rewards alongside individual immediate rewards, enhancing performance in scenarios where agents have limited access to partial information. However, their primary focus was on improving cooperation without considering training efficiency. Furthermore, Ackerman et al. [40] extended TD3 to multiagent scenarios to address value function overestimation. Mao et al. [41] introduced attention layers in critics, whereas Wang et al. [42] established communication channels using recurrent neural networks for partially observable environments. Iqbal et al. [43] utilized attention mechanisms to streamline information, thereby enhancing scalability. Yet, these approaches either increase computational complexity or struggle to handle the stochastic nature of the problem. Additionally, they treat all state information at the same importance level.

Consequently, this paper proposes an enhanced multiagent deep deterministic policy gradient (MADDPG) algorithm that tackles this concern by incorporating the system property and separately processing the crucial information. The critic network specifically accommodates the most important information i.e., stochastic market demand, as a distinct input alongside the other state inputs. This separation ensures that the most critical details receive special attention during the learning process, as they significantly impact the decision-making for efficient scheduling. By integrating the system property and processing this essential information separately, the algorithm focuses on prioritizing and accurately utilizing the most vital data, optimizing the scheduling strategy based on its significance in real-time decision-making scenarios.

## 3. System Description

In this paper, a multi-product FMS is considered, which consists of $M$ machines, $M-1$ intermediate buffers, $l$ product types, and $n$ mobile robots for loading and unloading parts. To illustrate this system, an example is provided in Fig. 1. The example showcases a configuration with 4-machines (denoted as $S_{ij}^l, l = 1,2,3 \ and \ i,j = 1,2,3,4$ represented by rectangles, where $i$ and $j$ represents product based sequence and the actual sequence respectively), 3-buffers (denoted as $b_i$, $i = 1,2,3$ represented by circles), 3-product types (denoted as $l = 1,2,3$ represented by small colored circles), 2-mobile robots (denoted as $R_i$, $i = 1,2$ represented as robot icons), a source for raw products, and a sink for finished products. In the system, after a machine completes its processing, a robot unloads the partially processed product to the corresponding downstream buffer based on the product type. Subsequently, the robot loads a new part from the upstream buffer. To demonstrate product types and their associated processing routes, color coding for the flow representation is used. The dashed rectangles represent the slowest machine (highest cycle time) for each product type. The flowlines linking stations and buffers, distinguished by different colors, illustrate the processing flows for the respective product types.

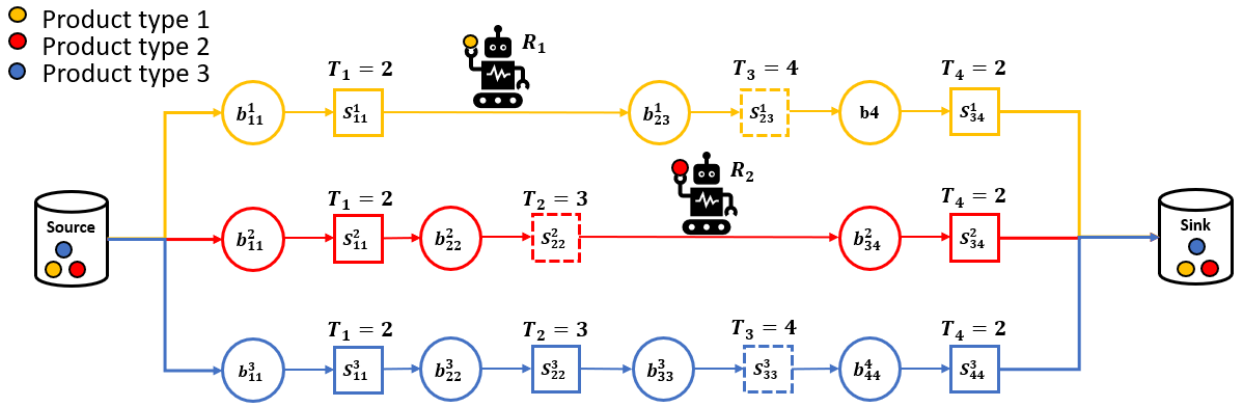

Fig.1. A multiproduct FMS with 4-machines, 3-buffers, 3-product types, and 2-robots

A summary of the notations used in the environment is presented in Table 1. The following constraints are considered in this study:

- The production line is capable of processing $l$ different product types, each with its unique sequence of operations.
- Every part must go through the first and last machine ($S_1$ and $S_M$).
- Every machine must process at least one product type.
- Part and product type are used interchangeably in this paper;
- The intermediate buffer only contains product types meant to be processed by the next machine.
- A machine is considered blocked if it is operational, but its downstream buffer is full, and considered starved if it is operational but its upstream buffer is empty.

**Table 1**
Table of Notations

| Symbol | Definition |
|---|---|
| $S_i$ | Machine $i$, where $i = 1, 2, \ldots, M$ |
| $S_{sc}$ | Slowest critical machine |
| $B_i$ | Buffer $i$ (Also denote buffer capacity), where $i = 1, 2, \ldots, M-1$ |
| $l$ | Product types, $l = 1, 2, \ldots, K$ |
| $\vec{b}_i(t)$ | Level of buffer $i$ at time $t$, $\vec{b}_i(t) = [p_{1i}(t) \ldots p_{li}(t) \ldots p_{Ki}(t)]'$ |
| $T_i$ | Processing time of $S_i$ |
| $T_{load}^i$ | Loading time on $S_i$ |
| $T_{unload}^i$ | Unloading time on $S_i$ |
| $T_{travel}$ | Travel time of the robot between $S_i$ and $S_{i+1}$ |
| $T_{travel}^{ij}$ | Travel time between $S_i$ and $S_j$, $T_{travel}^{ij} = |j - i| \times T_{travel}$, $(i, j = 1, 2, \ldots, M)$ and $i \neq j$ |
| $\vec{e}_i$ | $\vec{e}_i = (j, t_k, d_k)$, Disruption event $j$ happened at a time $t_k$ and lasting for a duration of $d_k$ |
| $\vec{e}_{vi}$ | Virtual disruption event |
| $d_i(t)$ | Duration of downtime on the machine $S_i$ at time $t$ |
| $D_l(t)$ | Market demand of product type $l$ at time $t$ |
| $R_i$ | Mobile robot $i$, where $i = 1, \ldots, K$ |
| $w_i(t)$ | Disturbance $i$ at time $t$ |
| $u_{ij}^l(t)$ | The control input of mobile robot $i$ for machine $S_j$ and product type $l$ at time $t$ |
| $OW_i$ | Opportunity window of the machine $S_i$ |
| $PPL_i$ | Permanent production loss of machine $S_i$ |
| $a_t^i$ | The action of agent $i$ at time $t$ |
| $r_t^i$ | Reward of agent $i$ at time $t$ |
| $s_t$ | Environment state at time $t$ |
| $\gamma$ | Discount factor |
| $\alpha$ | Learning rate |
| $\alpha_a$ | Actor learning rate |
| $\alpha_c$ | Critic learning rate |
| $\pi_t^i$ | Policy of agent $i$ at time $t$ |
| $o_t^i$ | Observation of agent $i$ at time $t$ |
| $g_i(t)$ | Mobile robot's assignment status to machine $i$ at time $t$ |
| $\theta_i(t)$ | Running status of machine $i$ at time $t$ |
| $\vec{P}_t$ | Vector of product types loaded onto the machines at time $t$ |
| $CT$ | Cost of demand dissatisfaction |
| $CoD(t)$ | Cost of delay at time $t$ |
| $\omega_d^l$ | Delay cost rate for product type $l$ |
| $CoS(t)$ | Cost of surfeit at time $t$ |
| $\omega_s^l$ | Surfeit cost rate for product type $l$ |
| $Y_l(t)$ | Last machine production (system's throughput) for product type $l$ at time $t$ |

### 3.1. Introduction of System Model and System Property

In our previous research work [26], a data-driven model was developed for a multiproduct FMS to effectively assess the current states of the system in real time. To provide a comprehensive understanding of the current paper, a concise overview of the main findings and the modelling approach without delving into the detailed derivation will be presented. The multiproduct FMS can be modelled with the following state-space equation:

$$\dot{b}(t) = F(b(t), U(t), W(t)) \tag{1}$$

In the context of this system, the parameters are defined as

$b(t) = [\vec{b}_1(t), \vec{b}_2(t), \cdots, \vec{b}_{M-1}(t)]'$ represents the buffer levels at time $t$.

$W(t) = [w_1(t), w_2(t), \cdots, w_M(t)]'$ represents the disturbances at time $t$, where $w_j(t)$ describes whether $S_j$ suffers from a disruption at time $t$. If $\exists \vec{e}_k \in E. s.t. \vec{e}_k = (i, t_k, d_k)$ and $t \in [t_k, t_k + d_k]$, then, $w_i(t) = 1$, otherwise, $w_i(t) = 0$. Define $\theta_i(t)$ as the status of a machine $S_i$ at time $t$ i.e., $\theta_i(t) = 1 - w_i(t)$. A machine $S_i$ is up at time $t$ when $w_i(t) = 0$, and down when $w_i(t) = 1$;

$$U(t) = \begin{bmatrix} u^1_{11}(t) & u^2_{12}(t) & \cdots & u^l_{1M}(t) \\ u^1_{21}(t) & u^2_{22}(t) & \cdots & u^l_{2M}(t) \\ \vdots & \vdots & \vdots & \vdots \\ u^1_{K1}(t) & u^2_{K2}(t) & \cdots & u^l_{KM}(t) \end{bmatrix}$$ is the control input, where $u^l_{ij}(t)$ describes whether the robot $i$ is assigned

to the machine $S_j$ for loading/unloading of product type $l$. There must be only one robot assigned to a machine at a time $t$. Although $U(t)$ is presented as a 2D matrix, it's crucial to note that it represents a 3D tensor with $i, j$, and $l$ axis representing the index for robots, machines and product types. Each element of this tensor is a binary variable $u_{ijl}(t)$, which is equal to 1 if a mobile robot $i$ is assigned to machine $j$ to load product type $l$ at time $t$, otherwise is 0.

The mobile robot's function comprises two major steps: unloading and loading the machine, as well as transporting the material. A mobile robot $i$ can only be assigned to machine $j$ if it is neither starved nor blocked and is not already assigned a mobile robot, i.e., the necessary condition for robot $i$ assignment is $(\sum_{j=2}^{M-2} b_j > 0) \& (B_{j+1} > \sum_{j=2}^{M-2} b_{j+1}) \& (u_{ij}(t-1) = 0)$.

A recursive algorithm has been developed, based on the principle of flow conservation, to evaluate the buffer level $b(t)$ at each time step $t$ in the system. The specific dynamic function $F(*)$ that describes this algorithm can be found in the referenced work [26]. To effectively control the assignments of robots in real time, it is essential to have a metric that assesses the system's dynamic performance. Prior research [26] has shown that not all downtimes result in permanent production loss (PPL) for the system. PPL only occurs when the downtime causes the slowest critical machine, denoted as $S_{sc}$, to stop. To quantify this characteristic, the concept of an opportunity window (OW) has been introduced. Each machine, denoted as $S_i$, has an associated OW, denoted as $OW_i$, which represents the maximum allowable downtime for $S_i$ without causing PPL for the system [26]. It has been established that $OW_i$ corresponds to the duration required for all buffers between $S_i$ and the slowest critical machine $S_{sc}$ to either become empty or full [26]. According to [26], the system experiences PPL only if any disruption event exceeds the OW of the slowest critical machine $S_{sc}$. Therefore, the occurrence of PPL during a given time interval $[0, T]$ depends on the status of $S_{sc}$ and can be defined as $PPL(T) = \frac{D_{sc}(T)}{T_r^{sc} + T_{load}^{sc} + T_{unload}^{sc} + T_{travel}^{j,sc}}$, where $D_{sc}(T)$ represents the duration of downtime for $S_{sc}$.

Detailed derivation and rigorous proof of OW and PPL can be found in [26]. $OW_i$ and $PPL_i$ serve as significant dynamic properties of the system and are utilized in designing the control policy presented in this paper.

### 4. Control Problem Formulation

Controlling production systems that involve multiple product types and stochastic market demand poses a greater challenge due to the increased complexity. In this research paper, the control actions entail assigning a robot, $R_i$ to a machine $S_j$ to load product type $l$. The absence of a closed-form representation for the system renders the application of traditional control methods challenging. Moreover, the problem's complexity, involving numerous states and unknown transition probabilities, hinders the use of operations research techniques like dynamic programming and centralized control schemes. Consequently, the problem is formulated as a Decentralized Partially Observable Markov Decision Process (Dec-POMDP) and addressed using Multi-Agent Reinforcement Learning (MARL). The control problem aims to determine an optimal policy in stochastic scenarios, mapping states to actions to maximize the reward.

Dec-POMDP is a framework employed in multi-agent environments for decision-making and coordination, represented by the tuple $(n, S, \{U_i\}, \{O_i\}, P, r, \gamma)$ [44] utilized in multi-agent environments for decision-making and coordination. In this framework, there are $n$ ageins denoted by $i \in \{1,2,...,n\}$ collaborating to achieve a cooperative task. At each decision point, agent $i$ selects an action $u_i \in U$ forming a joint action $\mathbf{u} \in \mathbf{U} \equiv U^i, \forall i \in \{1,2,...,n\}$. This action influences the true state of the environment $s \in S$ according to a transition probability function $P(s'|s, \mathbf{u}): S \times U \times S \to S$. Consequently, a global reward function $r(s, \mathbf{u}): S \times U \to \mathbb{R}$ is generated. Each agent has access only to partial observations $o^i \in O^i$ and selects actions based on a policy $\pi(\mu^i|\sigma^i)$ that depends on its local observation-action history $\sigma^i \in \Sigma \equiv (O^i \times U)$. The objective of each agent is to maximize the discounted sum of rewards over an episode, using a discount factor $\gamma$.

### 4.1. Dec-POMDP Framework for MARL

In MARL, control decisions typically involve generating a probability distribution over the action space. At each time step $t$, one action is sampled for each agent, resulting in a new state and a corresponding reward. The transition to the new state is governed by a transition probability determined by the system. Agents learn to select actions based on their observations of the environment, known as the policy. The quality of the policy is assessed based on the discounted sum of rewards obtained over an episode, and the objective of MARL is to compute a policy that maximizes this sum. Since the transition probabilities are unknown in this system, a model-free MARL algorithm is utilized. This approach eliminates the need for assumed transition probabilities and focuses solely on defining observations, actions, and rewards.

#### 4.1.1. Observations

The observation of each agent $i$ is represented as,

$$o_t^i = \{\vec{g}_t, \vec{P}_t\} \quad (2)$$

where $\vec{g}_t = [g_1(t), g_2(t) ..., g_M(t)]$ is the status of machines e.g., $\vec{g}_t = [0,1,0,0]$ represents that only a second machine is assigned a robot. $\vec{P}_t$ is a vector representing the product types each machine is currently processing e.g., in the case of 3-product types and 4-machines, $\vec{P}_t = [[0,1,0], [1,0,0], [0,0,0], [0,0,1]]$ represents that the first machine is loaded with product type 1, the second with product type 0, the third machine is empty, and the fourth one is loaded with product type 2.

In a Dec-POMDP, the state usually comprises the observations of all agents i.e.,

$$s_t = \{o_t^1, o_t^2, ..., o_t^n\} \quad (3)$$

However, in this case, besides the observations of all agents, additional information is provided. The state information $s_t$ at time $t$ can be presented as,

$$s_t = \{\vec{b}_t, \vec{\theta}_t, \vec{g}_t, \vec{P}_t\} \quad (4)$$

where $\vec{b}_t, \vec{\theta}_t, \vec{g}_t,$ and $\vec{P}_t$ are vectors representing buffers' levels, machines' status i.e., on/off, machines' assignment status i.e., assigned robot or not, and type of product loaded onto individual machines respectively.

#### 4.1.2. Action

The action an agent can take is the assignment action which defines the product type $l$, to be assigned to a machine $S_j$. It is defined as

$$A = \{a^i_{lj}\} \tag{5}$$

where $a^i_{lj}$ represents the action of agent $i$ to load product type $l$ to machine $j$ and $a^i_{lj} \in [0,1]$. A new action is triggered when a machine is empty or finishes processing a product and a robot is available.

*4.1.3. Reward*

A global reward characterizes the effectiveness of the collective action of all agents. Therefore, a suitable reward setting is crucial to ensure that the resulting policy aligns with expectations. Designing a good reward setting necessitates domain knowledge of the system and should be crafted to enhance the overall discounted sum of this reward, encouraging the system to function as intended. Various potential reward settings were investigated for the problem, and the following was determined to be the most effective.

$$R = -PPL(t) - CoD(t) + CoS(t) \tag{6}$$

where the first term punishes the system for accruing PPL, the second term represents the cost of delay, and the third term is the cost of surfeit i.e.,

$$PPL(T) = \frac{D_{sc}(T)}{T_r^{sc} + T_{load}^{sc} + T_{unload}^{sc} + T_{travel}^{j,sc}} \tag{7}$$

$$CoD(t) = \sum_{l=1}^{K} \omega_d^l \max\{0, D_l(t) - Y_l(t)\} \tag{8}$$

$$CoS(t) = \sum_{l=1}^{K} \omega_s^l \max\{0, Y_l(t) - D_l(t)\} \tag{9}$$

The agent is punished when the production is less than the corresponding demand at time $t$. Similarly, the agent is rewarded when the demand is satisfied. $\omega_d^l$ and $\omega_s^l$ are the delay cost rate and the surfeit cost rate for product type $l$, respectively.

*4.2. MARL-Based Control of Multiproduct FMS*

In this section, a decentralized control is introduced for a multiproduct FMS using a cooperative MARL framework. In this approach, each mobile robot is considered as an independent agent. These agents acquire the ability to interpret their individual observations, transforming them into assignment actions that reduce the occurrence of PPL while fulfilling the unique customer demands associated with each product type.

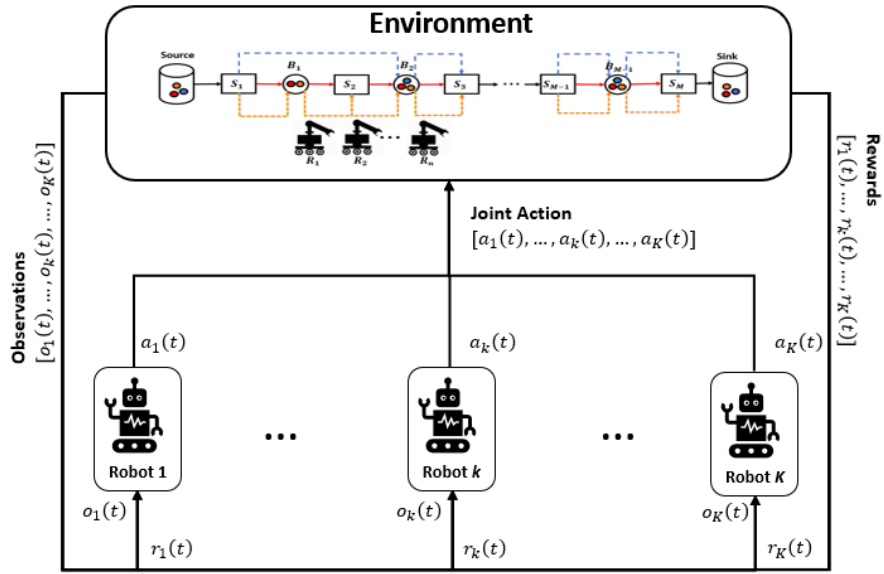

Fig.2. MARL framework for a multiproduct FMS

The mobile robots, acting as agents, sequentially interact with the multiproduct FMS within discrete-time steps to learn optimal assignment policies. The agent $k \in \mathcal{K}$ possesses its specific local observation space $\mathcal{O}_k$ and the overall environment state space is denoted as $\mathcal{S}$. Illustrated in Fig.2, at each time step $t$, the agent utilizes its local observation $o_k(t): \mathcal{S} \mapsto \mathcal{O}_k$ to determine an action $a_k(t) \in \mathcal{A}_k$, sourced from its action space $\mathcal{A}_k$, based on the current policy $\mu_k$. The shared environment receives the collective action of all agents $a(t) = (a_1(t), a_2(t), \dots, a_K(t))$ and transitions to the subsequent set of observations $o_k(t+1) \in \mathcal{O}_k$ while providing rewards $r_k(t): \mathcal{S} \times \mathcal{A}_k \mapsto \mathbb{R}$ for all $k \in \mathcal{K}$. The agents aim to enhance their individual policies continuously until convergence to the optimal policy $\mu_k^*$ that maximizes the expected long-term discounted cumulative reward, defined as $J_k(\mu_k) = E[\sum_{t=1}^{T} \gamma^{t-1} r_k(t)]$, where $\gamma \in [0,1]$ is the discount factor and $T$ represents the total number of time steps. The optimal policy $\mu_k^*$ can be obtained as

$$\mu_k^* = \arg\max_{\mu_k} J_k(\mu_k) \tag{10}$$

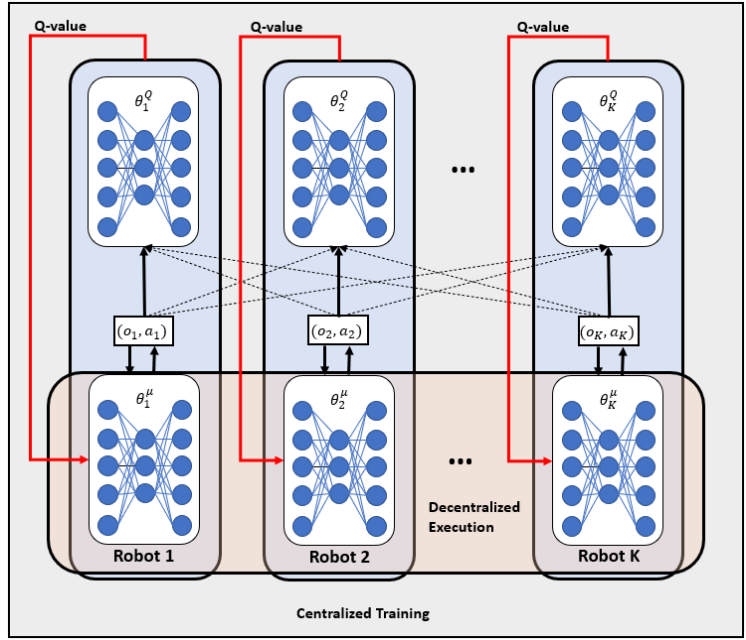

Fig.3. Centralized training and decentralized execution (CTDE) in MADDPG algorithm [45]

However, in a multi-agent scenario, the environment is no longer stationary from the viewpoint of an agent, given that other agents are concurrently updating their policies. This presents a moving target problem for the agents, leading to potential learning instability. To tackle these challenges, researchers have introduced the MADDPG method [45].

MADDPG is an actor-critic algorithm that utilizes policy gradients. During the training phase, multiple agents undergo centralized training, allowing for the exchange of additional information to enhance the training process. Illustrated in Fig. 3, the agents not only receive their local observations but also supplementary information, encompassing the observations and actions of other agents $a(t) = (a_1(t), a_2(t), \dots, a_K(t))$, and complete observation of the environment state $s_k(t) = (o_k(t), o_{-k}(t))$, where $o_{-k}(t)$ denotes the local observations of other agents at time step $t$. This extra information, including the actions taken by other agents, aids the agents in addressing the challenges posed by the dynamic nature of the environment, enabling them to comprehend and model the dynamics of the multiproduct FMS effectively. It is noteworthy that this additional information is not utilized during the execution phase, where agents solely rely on their local observations, ensuring a fully distributed decision-making process.

Agent $k$ in MADDPG utilizes two primary deep neural networks: an actor-network with parameters $\theta_k^\mu$ to approximate a joint policy $\mu_k(o_k|\theta_k^\mu)$ for robot scheduling and a critical network with parameters $\theta_k^Q$ to approximate the state-value function $Q_k(s_k, a|\theta_k^Q)$, but this network doesn't take into account the system property and cannot differentiate between the most important information i.e., changing market demand and less important information i.e., machines' and buffers' status. Treating it equally and processing it in the same layers results in longer training times. Therefore, a novel critic-network architecture is set up with a parallel layer for the varying market demand, as shown in Fig.4. By incorporating a parallel layer, the critic network can effectively process the stochastic market demand separately from the rest of the state information. This allows the network to focus on the specific characteristics and dynamics of market demand, leading to improved understanding and modelling of its impact on the system. These changes enable the critic network to capture the complexities and uncertainties associated with market fluctuations. This, in turn, facilitates better decision-making by mobile robots, as they can incorporate up-to-date information into their actions and policies. With this update in critic, the state-value function is now approximated with demand information $Q_k(s_k, d, a|\theta_k^Q)$.

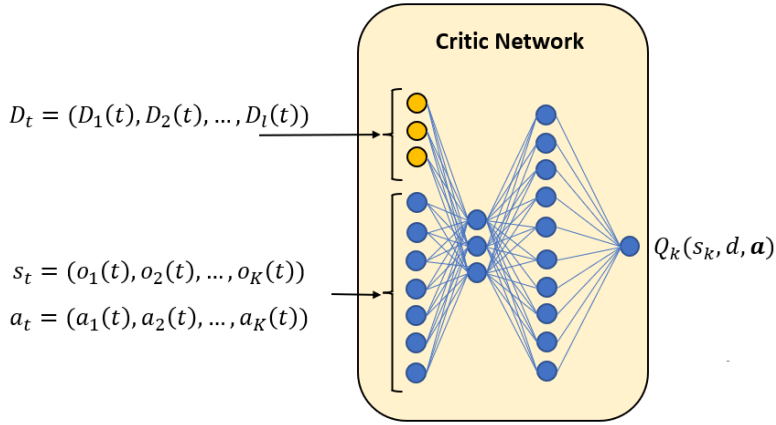

Fig.4. Modified Critic network of MADDPG algorithm

These networks also have corresponding time-delayed copies, $\theta_k^{\mu'}$ and $\theta_k^{Q'}$, which serve as target networks. At each time step $t$, the actor determines the continuous action $\mu_k(o_k(t)|\theta_k^\mu)$ based on the local observation $o_k(t)$, incorporating a random noise process $\mathcal{N}_k$ to introduce exploration. This results in the action $a_k(t) = \mu_k(o_k(t)|\theta_k^\mu) + \mathcal{N}_k(t)$. The shared environment gathers the joint action $a(t) = \{a_k(t), \forall k \in \mathcal{K}\}$, offering immediate reward $r_k(t)$, and subsequent observation $o_k(t+1)$ for each agent. To enhance sample efficiency and training stability, the agent stores its transition, along with additional information $e_k(t) = (s_k(t), a(t), r_k(t), s_k(t+1))$, in its replay buffer $\mathcal{D}_k$ for future decision-making steps.

To train the primary networks, a random mini-batch of $B$ samples $\left(s_k^i, a^i, r_k^i, s_k^{i+1}\big|_{i=1}^B\right)$ is selected from the replay buffer $\mathcal{D}_k$. The critic network undergoes an update by minimizing the loss function:

$$\mathcal{L}_k(\theta_k^Q) = \frac{1}{B}\sum_i \left(y_k^i - Q_k(s_k^i, d, a^i|\theta_k^Q)\right)^2 \tag{11}$$

where $y_k^i$ is the target value expressed as $y_k^i = r_k^i + \gamma Q_k'\left(s_k^{i+1}, d, a_{i+1}\big|\theta_k^{Q'}\right)\big|_{a_{i+1}=\{\mu_k'(o_k^{i+1}), \forall k \in \mathcal{K}\}}$, with a discount factor $\gamma$ and stochastic market demand $d$. The critic network's parameters, $\theta_k^Q$, are adjusted using gradient descent of Eq. (12) such that $\theta_k^Q \leftarrow \theta_k^Q - \beta_Q \nabla_{\theta_k^Q} \mathcal{L}_k(\theta_k^Q)$, with a learning rate of $\beta_Q$. Likewise, the policy network adjusts its parameters to maximize the expected long-term discounted reward $J(\mu_k|\theta_k^\mu)$ using the deterministic policy gradient, $\nabla_{\theta_k^\mu} J(\mu_k|\theta_k^\mu)$, such that $\theta_k^\mu \leftarrow \theta_k^\mu - \beta_\mu \nabla_{\theta_k^\mu} J(\mu_k|\theta_k^\mu)$, where

$$\nabla_{\theta_k^\mu} J(\mu_k|\theta_k^\mu) \approx \frac{1}{B}\left[\sum_i \nabla_{a_k} Q_k(s_k^i, d, a_i|\theta_k^Q) \times \nabla_{\theta_k^\mu} \mu_k(o_k^i|\theta_k^\mu)\right]\bigg|_{a_i=\{\mu_k(o_k^i), \forall k \in \mathcal{K}\}} \tag{12}$$

and $\beta_\mu$ denotes the actor-network's learning rate. The target parameters in both networks undergo updates through soft updates at a rate of $\tau$, as outlined below:

$$\theta_k^{\mu'} \leftarrow \tau\theta_k^\mu + (1-\tau)\theta_k^{\mu'} \tag{13}$$

$$\theta_k^{Q'} \leftarrow \tau\theta_k^Q + (1-\tau)\theta_k^{Q'} \tag{14}$$

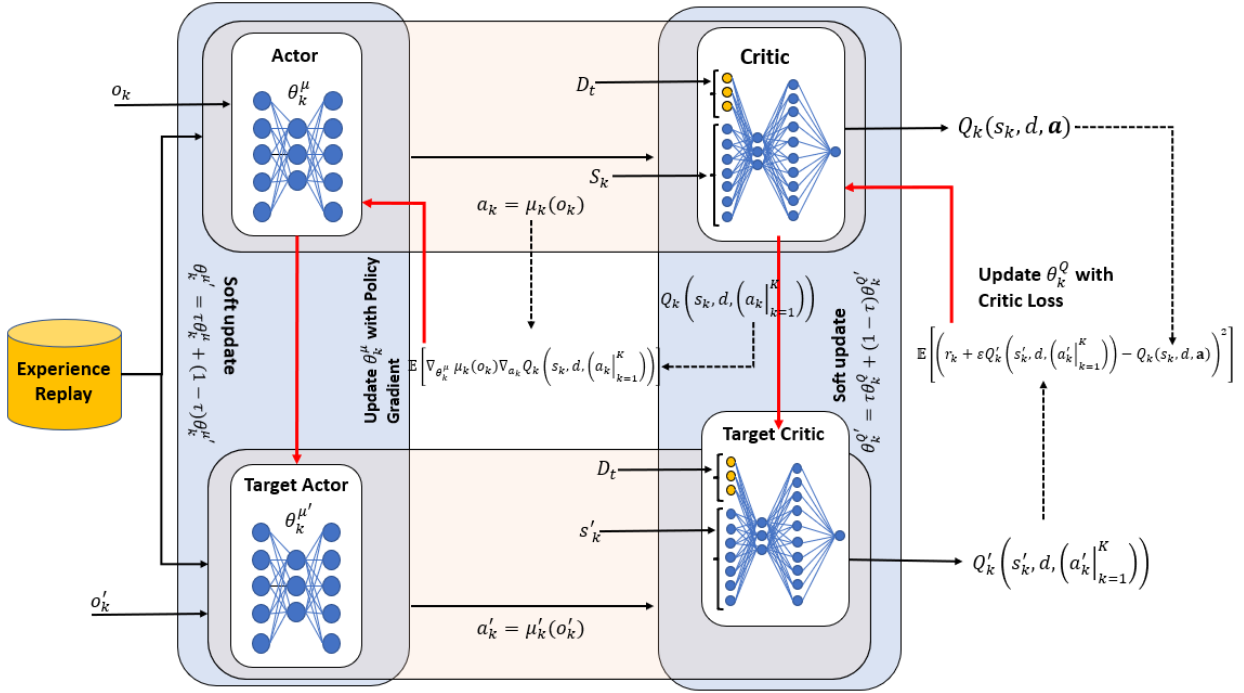

Fig.5. Representation of the training procedure in the modified-MADDPG algorithm for the $kth$ agent

where $\tau$ is a constant close to zero.

Fig. 5 illustrates the interactions among the four networks during the training process of the $k-th$ agent in the modified- MADDPG.

**Algorithm 1: Modified-MADDPG Algorithm**

**for** each agent $k \in \mathcal{K}$ **do**
    Initialize replay buffer $\mathcal{D}_k$
    Initialize the actor-network $\mu_k(o_k|\theta_k^\mu)$ and critic network $Q_k(s_k, d, a|\theta_k^Q)$ with weights $\theta_k^\mu$ and $\theta_k^Q$, respectively
    Initialize the target networks, $\mu'_k$ and $Q'_k$ with weights $\theta_k^\mu$ and $\theta_k^Q$, respectively
**for** each episode $e = 1, 2, \ldots$ **do**
    **for** each agent $k \in \mathcal{K}$ **do**
        Initialize random process $\mathcal{N}_k$ for exploration
        Generate initial local observation from the environment simulator
    **for** each step $t = 1, 2, \ldots$ **do**
        **for** each agent $k \in \mathcal{K}$ **do**
            Select action $a_k(t) = \mu_k(o_k(t)|\theta_k^\mu) + \mathcal{N}_k(t)$
        Execute joint action $a(t) = (a_1(t), a_2(t), \ldots, a_K(t))$
        **for** each agent $k \in \mathcal{K}$ **do**
            Collect reward $r_k(t)$ and observe $s_k(t+1)$
            Store the transition $(s_k(t), a(t), r_k(t), s_k(t+1))$ into $\mathcal{D}_k$
            Sample random minibatch of B transitions $(s_k^i, a^i, r_k^i, s_k^{i+1})$ from $\mathcal{D}_k$

Update the critic network by minimizing the loss given by (19)
Update the actor policy using the sampled policy gradient given by (20).
Update the target networks according to (21)

The stepwise implementation of modified-MADDPG is given in Algorithm 1.

## 5. Case Study

In this section, simulation experiments are conducted to validate the effectiveness of the proposed modified MADDPG algorithm. The experiments involve comparing the production output by using control policies learned from DQN, DDPG, simple MADDPG, and our modified MADDPG algorithm. DQN and DDPG are single-agent RL algorithms where a single agent controls all the assignments of mobile robots at each state, while MADDPG utilizes multiple agents. The algorithms are evaluated based on two performance metrics: (1) the training time required in obtaining a scheduling policy for the robots, and (2) the system throughput achieved using the trained policy for robot assignments. We aim to demonstrate two key findings: (1) with the same training duration, our proposed method achieves the highest performance in improving system throughput and demand satisfaction, and (2) our proposed method outperforms other control algorithms in terms of training process efficiency, reaching a stable reward faster.

### 5.1. Experiment Parameter Setting

The simulation experiments are conducted using the Python platform. We create and analyze 100 distinct production lines, which encompass various types such as automotive assembly lines, semiconductor lines, and battery lines, among others. Through these experiments, we can verify the effectiveness and resilience of the proposed method across 100 diverse production lines. Each line is constructed by randomly selecting from the following sets:

$$No.\ of\ machines, M \in \{3,4,\ldots,10\}$$

$$No.\ of\ robots, R \in \{1,2,3,4,5,6\}$$

$$Processing\ time, T_i \in [50, 90]\ sec, i = 1,2,\ldots,M$$

$$Loading\ time, T^i_{load} \in [5, 25]\ sec, i = 1,2,\ldots,M$$

$$Unloading\ time, T^i_{unload} \in [5, 25]\ sec, i = 1,2,\ldots,M$$

$$MTBF_i \in [10, 25]\ mins, i = 1,2,\ldots,M$$

$$MTTR_i \in [2, 6]\ mins, i = 1,2,\ldots,M$$

$$Buffer\ capacity, B_i \in [10, 50], i = 1,2,\ldots,M-1$$

$$No.\ of\ product\ types, l \in \{2,3,4,5,6\}$$

To illustrate the application of the proposed method, a comprehensive numerical case study involving a multiproduct flexible production line is provided. This production line consists of four machines, two robots, three product types, and three buffers, similar to the configuration shown in Fig. 1. The specific parameters for this setup are presented in Tables 2 and 3 below.

**Table 2**

Parameters for the machines

| Parameters | S1 | S2 | S3 | S4 |
|---|---|---|---|---|
| Processing time, Tr (sec) | 60 | 55 | 70 | 65 |
| Loading time, Tl (sec) | 15 | 15 | 15 | 15 |
| Unloading time, Tu (sec) | 10 | 10 | 10 | 10 |
| MTBF (min) | 16 | 14 | 22 | 20 |
| MTTR (min) | 5 | 3 | 4 | 4 |

Table 3
Parameters for the three buffers

| Parameters | B1 | B2 | B3 |
|---|---|---|---|
| Buffer Capacity | 25 | 25 | 25 |
| Initial Buffer Level | 0 | 0 | 0 |

Using the specified system parameters, Algorithm-1 is employed to train the control policy. The control decision is activated whenever a machine completes processing and is ready for unloading/loading. The neural network architecture consists of two fully connected hidden layers, each containing 64 hidden units. The reward function is formulated according to Eq. (6). The discount factor $\gamma$ for model value expansion is set to 0.95. The actor learning rate $\alpha_a$ is assigned a value of 0.0005, while the critic learning rate $\alpha_c$ is set to 0.001. Notably, the first layer of the critic network is accompanied by an additional parallel layer, which takes stochastic market demand as input. Upon completion of the training process, a machine control policy is derived based on the updated network parameter $\theta$.

To assess the performance of training and the learned policy $\pi$ obtained using the proposed method, it is compared with three other control strategies: DQN, DDPG, and simple MADDPG.

*DQN:* It combines Q-learning with deep neural networks to handle high-dimensional state spaces. It uses an experience replay buffer and a target network to stabilize and improve the learning process.

*DDPG:* It is an actor-critic algorithm designed for continuous action spaces. It employs an actor-network to approximate the policy and a critic network to estimate the value function. It uses a replay buffer and target networks to enhance stability during training.

*MADDPG:* It extends DDPG to multi-agent environments, allowing multiple agents to learn concurrently by incorporating their observations and actions into the training process. MADDPG addresses the challenges posed by non-stationarity and coordination among agents by leveraging centralized training with decentralized execution.

6. Results and Discussion

*6.1. Training performance comparison*

The training processes of all control algorithms are depicted in Fig. 6. The x-axis represents the training episodes, while the y-axis corresponds to the rewards defined in Eq. (6). It's noticeable that both the simple MADDPG and the proposed modified-MADDPG algorithms eventually attain higher rewards in comparison to the DDPG and DQN algorithms. Overall, the DQN algorithm necessitates a relatively longer training period to stabilize and yield the lowest reward. Conversely, the DDPG algorithm swiftly adapts to the control task, achieving stability, although with a lower reward than the simple and modified MADDPG algorithms. Particularly, the modified-MADDPG algorithm demonstrates the highest reward and achieves stability rapidly. The simple MADDPG algorithm takes slightly longer than DDPG and modified-MADDPG to stabilize, yet it outperforms DDPG in terms of reward. These results validate that the modified MADDPG algorithm exhibits superior performance in acquiring the optimal control policy.

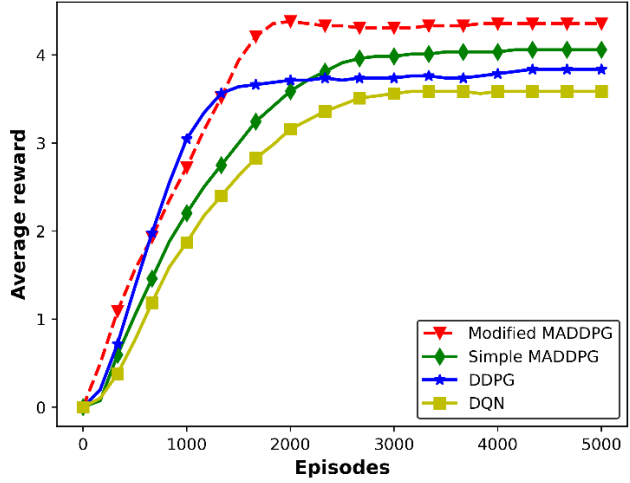

Fig.6. Training process of multiproduct FMS with DQN, DDPG, simple MADDPG, and modified-MADDPG algorithms

### 6.2. Execution performance comparison

Figures 7a and 7b present a comparison of the average production levels, accompanied by their respective 95% confidence intervals, achieved by the trained policies using the four control algorithms. The control policy is significantly influenced by market demand, playing a crucial role in determining production outcomes. The analysis covers two scenarios: Fig. 7a assumes an equal probability for the arrival of demand for each product type, while Fig. 7b assigns higher weightage to product type 2 in contrast to product types 1 and 3. In both scenarios, it becomes evident that the proposed modified MADDPG algorithm surpasses the other methods in effectively meeting customer demand.

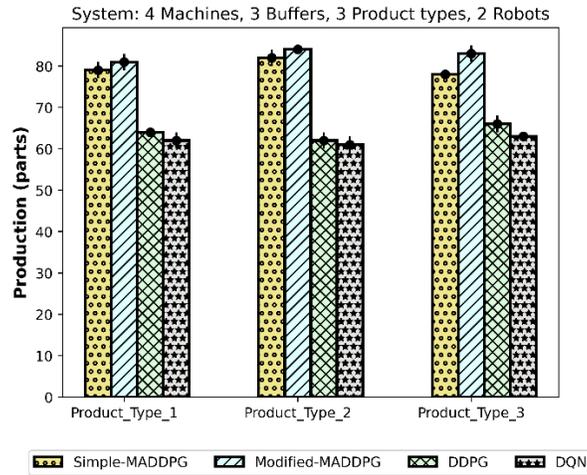

Fig. 7(a). Production comparison under the four control algorithms with equiprobable demand of individual product type

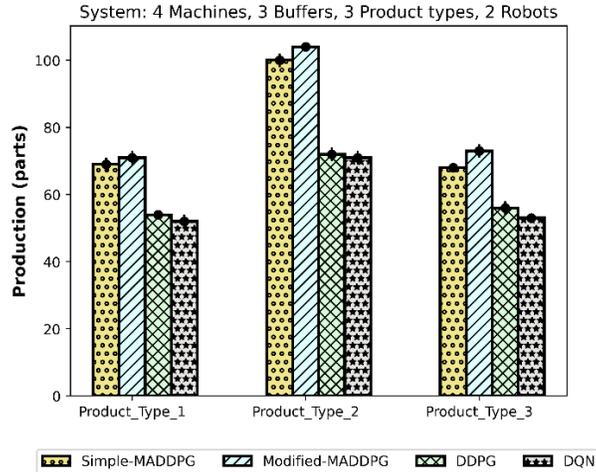

Fig.7(b). Production comparison under the four control algorithms with un-equiprobable demand of individual product type

In Figure 7a, the proposed method achieves an average production of 85 units for each product type, outperforming the simple MADDPG algorithm which attains an average of 80 units for each type. Conversely, the DDPG and DQN methods yield lower production levels, averaging around 62 units per product type. Shifting to Figure 7b, all control algorithms adeptly adapt their production levels based on market demand, assigning a larger quantity to product type 2 compared to the other types. Nevertheless, the comparative analysis remains consistent: the proposed method displays higher production levels than the other methods, with the simple MADDPG algorithm also surpassing DDPG and DQN. Although the production difference between DQN and DDPG is not substantial, DDPG shows slightly better performance.

## 7. Transfer Learning for Complex Environments

The proposed method, in conjunction with other control algorithms, exhibits notable results in the specific environment configuration, as depicted in Figure 1. However, exploring the effectiveness of this proposed approach in more complex environments remains crucial, given its potential variability in efficacy. Consequently, further investigations have been undertaken to assess the adaptability of the proposed method in complex environments.

In a broader scenario, a complex environment can be randomly generated, encompassing a larger array of machines, buffers, robots, and product types. Training the robot scheduling policy for such a complex environment poses challenges, particularly regarding the time-intensive nature of training the proposed algorithm from scratch. To mitigate this challenge, one potential solution involves leveraging a pre-trained model from a similar environment.

To exemplify this approach, a multiproduct FMS featuring 8 machines, 7 buffers, 3 product types, and 4 mobile robots serves as a new complex environment. Utilizing similarities shared with the configuration shown in Fig. 1, the previously trained model is applied to train the robot scheduling policy for this new complex environment. The training process for this updated configuration is illustrated in Fig. 8, comparing the training outcomes with and without incorporating the trained layers from the prior environment configuration.

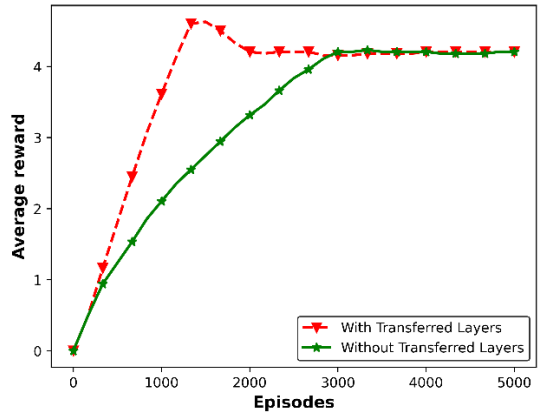

Fig.8. Training of proposed MADDPG algorithm with and without pre-trained layers in a complex environment

The incorporation of a trained model from a previous environment configuration noticeably expedites the training process, facilitating the establishment of a stable robot scheduling policy in more complex environments. This underscores the potential efficacy of reusing trained models to accelerate training and attain effective control policies within intricate environments.

## 8. Conclusion and Future Work Recommendations

This research paper focuses on addressing the scheduling policy for multi-agent mobile robots in a multiproduct Flexible Manufacturing System (FMS) environment, which encounters random machine disruptions and stochastic market demand. The problem of scheduling robots is formulated as a Multi-Agent Reinforcement Learning (MARL) problem, and a decentralized Modified Multi-Agent Deep Deterministic Policy Gradient (MADDPG) method is proposed as a solution to obtain an optimal robot schedule.

The proposed method introduces modifications to the MADDPG algorithm by incorporating a parallel layer in the critic network, specifically designed for efficient processing of stochastic market demand. The performance evaluation of the system is based on metrics such as market demand satisfaction and system throughput. The advantages of modifying the MADDPG algorithm are thoroughly discussed. Comprehensive case studies are conducted in the context of a multiproduct FMS to demonstrate the effectiveness of the proposed method in achieving an optimal scheduling policy for robots in the presence of stochastic market demands. Furthermore, the generalizability of the method to complex environment configurations is established by showcasing its superiority over other machine learning algorithms such as simple MADDPG, DDPG, and DQN methods in terms of optimizing robot scheduling.

It's essential to note that this work assumes fixed cycle times for the machines and overlooks energy consumption. In real-world scenarios, machines need flexibility to operate at various speeds to adapt to market demands and optimize energy use. Therefore, future research will address the inclusion of variable cycle times, allowing dynamic updates based on market demand and energy consumption considerations.

**Declaration of Competing Interest**
    No potential conflict of interest was reported by the author(s).

**Funding**
    This work was funded by the U.S. National Science Foundation (NSF) Grants 1853454 and 2243930.